\RequirePackage{fix-cm}
\documentclass[smallcondensed]{svjour3}     
\smartqed  
\usepackage{graphicx}
\newcommand{\bT}{{\bf T}}
\newcommand{\bU}{{\bf U}} \newcommand{\bC}{{\bf C}}

\newcommand{\bP}{{\bf P}}

\newcommand{\bk}{{\bf K }}
 \newcommand{\bB}{{\bf B}}

\newcommand{\bH}{{\bf H}}
\newcommand{\bF}{{\bf F}}
\newcommand{\bQ}{{\bf Q}}
\newcommand{\bZ}{{\bf Z}}
 \usepackage{mathptmx}      
 \usepackage{amssymb,amsmath}
\usepackage{subfigure}
\usepackage{latexsym}
\begin{document}

\title{Direct numerical scheme for all classes of nonlinear Volterra integral equations of the first kind
}

\titlerunning{Direct operational vector scheme for NVIE1}

\author{  R. Dehbozorgi \and
      K. Maleknejad 
}


\institute{R. Dehbozorgi$^ * $ \and K. Maleknejad
        \at
             School of Mathematics, Iran University of Science and Technology, Narmak, Tehran 16844, Iran  \\
              Tel.: +98 21 732 254 16\\
              Fax: +98 21 730 216 62\\
              \email{R.dehbozorgi2012@gmail.com; Maleknejad@iust.ac.ir; }
}

\date{Received: date / Accepted: date}

\maketitle

\begin{abstract}
This paper presents a direct numerical scheme  to approximate the solution of all classes of nonlinear Volterra integral equations of the first kind. This computational method is based on operational matrices and vectors. The operational vector for  hybrid block pulse functions and Chebyshev polynomials is constructed.  The scheme transforms the integral equation to a matrix equation and solves it with a careful estimate of the error involved. The main characteristic of the scheme is the low cost of setting up the equations without using any projection method which is the consequence of using operational vectors. Simple structure to implement, low computational cost and perfect approximate solutions  are the major points of the presented method.
 Error analysis and comparisons with other existing schemes demonstrate the efficiency and the superiority of our scheme.
\keywords{Direct method \and Nonlinear Volterra integral equations of the first kind \and Operational vector \and Hybrid block pulse and Chebyshev polynomials.}
\subclass{MSC 65R20 \and MSC 45D05 \and 45G10}
\end{abstract}

\section{Introduction}
\label{intro}
 In literature, many numerical methods for solving integral equations of the second kind have been presented by several authors (Atkinson 1997, Brunner 2004, Aziza and Islam 2013, Conte and Paternoster 2009,  Sahu and Saha Ray 2014,  Ghoreishi and Hadizadeh 2009, Guru Sekar and Murugesan 2016, Maleknejad and Dehbozorgi 2018, Maleknejad et al. 2007).
   In comparison with the abundant research concerned
with the numerical analysis of these equations, a few computational approaches have
been established to approximate the solution of integral equations of the first kind, especially in the
nonlinear case. Nonlinear Volterra integral equations of the first kind (NVIE1) appears as a famous  mathematical model in physics and engineering problems, e.g. electrochemical systems (Bieniasz 2015), electrostatic (Ding et al. 2003), heat conduction
problems (Bartoshevich 1975), etc.
Regarding  this fact that Volterra integral equations of the first kind are inherently
ill-posed (slight changes in inputs make large errors in outputs), so choosing the best numerical schemes due to overcoming the difficulty of ill-posedness  is significant. Tikhanov and Arsenin (1977) proposed various regularization techniques to conquer ill-posedness but practically, obtaining an appropriate filter to regularize is too difficult and time-consuming.\\
However, several methods have been developed to solve these types of equations, but a few numerical methods can be conducted in the nonlinear case. For instance, the approaches presented by Masouri et al. (2010) and Babolian et al. (2008), solved the VIE1 by using the expansion-iterative method and  operational matrix method. Maleknejad et al. (2007, 2011) have presented numerical techniques
based on wavelets, modified block pulse functions and Bernstein’s
approximation method for the solution of VIE1, respectively. Khan et al. (2014) have described optimal Homotopy asymptotic
method for solving these equations.   As a challenge of overcoming the ill-posedness and the nonlinearity of these equations together,  Babolian and Masouri (2008) and Babolian et al. (2008)  have presented a direct method to solve some particular NVIE1 using the operational matrix with block-pulse functions and operational matrices of piecewise constant orthogonal functions, respectively.  Singh and Kummar (2015) have described Haar wavelet operational matrix method for a class of NVIE1. As the matter of fact, it seems that operational matrices play preconditioner role in this kind of equations, for more details see Masouri et al. (2010), Babolian and Masouri (2008), Maleknejad et al.(2007), Babolian et al. (2011), Sing and Kummar (2016). Recently the Sinc Nystr{\"o}m method  has been applied to solve these equations in Ma et al. (2016). \\
In the present paper, we introduce a direct computational method to determine the approximate solution of some classes of nonlinear Volterra integral equations of the first kind. This scheme consists of reducing these equations to a nonlinear system of algebraic equations by expanding the given functions as Chebyshev polynomials (CP) with unknown coefficients. In some cases for improving the accuracy, we utilize the hybrid block-pulse functions and Chebyshev polynomials (HCP), especially for the unknown solutions belong to the class $ C^0/C^1$. Regarding our previous work (Maleknejad and Dehbozorgi 2018), we provide the operational vector for HCP. This vector together with the operational matrix of integration and product are then utilized to evaluate the unknown coefficients. Operational vector eliminates one of  the basis functions vectors which yields to have a direct method instead of using any projection methods, specially collocation method which is unable to work well for   integral equations of the first kind, for instance, see  Conte and Paternoster (2009) and Maleknejad et al. (2007).
\\
The objective of this study is to propose an efficient numerical scheme based on operational vectors for solving the
following class of NVIE1:
 \begin{equation}\label{V}
f(t)=\int_{t_0}^t K(x,t)G(u(x)) dx, ~~~~t\in D:=[t_0,t_f]
\end{equation}
where $ K, G $ and $ f $ are the given smooth functions and $G$ is a nonlinear function in terms of the unknown function $u(x)$. It is assumed that $ f(t_0)=0. $ Some important forms of function $ G(u(t)) $ are as follows:
\begin{itemize}
\item$u^{(n)}(t)$, it is assumed that $u = u^{(1)} = u^{(2)} = ........ = u^{(n-1)} = 0$ at $t = t_0$, where
$ u^{(n)}$ represents the nth derivative of $  u$ with respect to $x $,
\item $ u^{\alpha}(t) $, the $ \alpha $th power of $u(t)$, $ \alpha \in \mathbb{R} $,
\item $sin(u(t)), cos(u(t)), ln(u(t))$ and $ e^{u(t)} $ or any combination of these functions,
\item $ \sum\limits_{r=0}^{m} \alpha_ru^r(t) , ~~m\in \mathbb{N}$.
\end{itemize}
These mentioned forms of $ G $ may classify as  invertible or algebraic nonlinear functions over some especial intervals. In this study, we propose different techniques for all forms of nonlinear function $ G(u(t)) $ .
\\The paper is organized as follows: In Section 2, we briefly state some basic concepts of CP and HCP. As the key idea, we introduce the  operational vectors for these polynomials. In Section 3, the outline of the scheme  is presented.  Some theorems for the error analysis are presented in Section 4. In Section 5, numerical results   verify the applicability of our method in comparison with other existing methods (Sing and Kummar 2016, Ma et al. 2016).
\section{Preliminaries}
\subsection{Chebyshev polynomials (CP)}
Chebyshev polynomials (CP) of the first kind are defined by (Abramowitz and Stegun 1970) as
$$\phi_{n}(x)=cos(n\theta),~~ \theta=Arc~cos(x),$$
The orthogonality condition for CP with the weight function $w(x)=(\sqrt{1-x^2})^{-1}$ is as follows
\begin{equation*}
\langle \phi_{i}(x),\phi_{j}(x)\rangle_{w}=\int_{-1}^{1}w(x) \phi_{i}(x)\phi_{j}(x) dx=\delta_{ij}\left\lbrace \begin{array}{cc}\vspace*{0.15in}
\dfrac{\pi}{2},~~~~i\neq 0,\\
\pi ,~~~~~i=0.
\end{array}\right.
\end{equation*}
where $ \delta $ is Kronecker delta.
\\Shifted Chebyshev polynomials (SCP) of degree $ m $  is defined over the interval $D=[t_0,t_{f}]$ as follows
\[T_{m}(t)=\phi_{m}(\dfrac{2}{t_{f}-t_{0}}(t-t_0)-1)=\phi_{m}(A(t-t_0)-1),\]
where
\begin{equation}\label{A}
A=\frac{2}{t_f-t_0},
\end{equation}
 and so the weight function for SCP is described as (Datta and Mohan 1995, p. 90)
\[\tilde{w}(t):=w(A(t-t_0)-1)=\dfrac{t_f-t_0}{2\sqrt{(t-t_0)(t_f-t)}}.\]
One of the important properties of CP is completeness, therefore SCP also form a complete orthogonal set, that is,  every $f \in L^2(D)$  can be represented as an infinite series $$f(t)=\sum \limits_{r=0}^{\infty} c_{r} T_{r}(t),$$ where the coefficient $c_{i}$ can be determined as
\begin{equation}\label{cr}
\hspace*{0.2in}~~~~~~~~c_{i}=\dfrac{\langle f(t),T_{i}(t)\rangle_{\tilde{w}}}{\langle T_{i}(t),T_{i}(t) \rangle_{\tilde{w}}}.
\end{equation}
Moreover, the orthogonality condition of SCP is as follows
\begin{equation}\label{wt}
\langle T_{i}(x),T_{j}(x)\rangle_{\tilde{w}}=\int_{t_0}^{t_f}\tilde{w}(x) T_{i}(x)T_{j}(x) dx=\frac{\delta_{i,j}}{A}\left\lbrace \begin{array}{cc}\vspace*{0.15in}
\pi ,~~~~i=0,\\\vspace*{0.15in}
\dfrac{\pi}{2},~~~i\neq 0.\\
\end{array}\right.
\end{equation}
where $ A $ is defined in Eq. (\ref{A}).
\subsection{Hybrid Chebyshev polynomials and block pulse functions (HCP)}
Hybrid Chebyshev polynomials and block pulse functions (HCP) $ T_{im}(t)$ have two parameters where $ n=1,...,N$ and $  m=0,1,...,M-1$ are
the order of block-pulse functions and Chebyshev polynomials, respectively. They are defined  over the interval $ D$ as
\begin{equation*}
T_{im}(t)=\left \lbrace \begin{array}{ll}
T_{_m}(AN(t-t_0)-2n+1),&~~~~~t\in \left[t_0+\frac{2(n-1)}{AN},t_0+\frac{2n}{AN}\right],\\
0,& ~~~~~otherwise.
\end{array}\right.
\end{equation*}
where $ A $ is defined in (\ref{A}) and $ T_{_m}(t) $, $m=0,1,...,M-1$ are the Chebyshev polynomials
which are defined over the interval $[-1,1]$. The weight functions for hybrid Chebyshev polynomials and block pulse function are
\begin{equation}
\hspace*{0.2in}\omega_n(t)=w((AN(t-t_0)-2n+1),\hspace*{0.5in}t\in [t_0+\frac{2(n-1)}{AN},t_0+\frac{2n}{AN}]
\end{equation}
where  $ w(t)=(\sqrt{1-t^2})^{-1}. $ In general case, $ \tilde{w}(t):=\sum\limits_{n=1}^N \omega_n(t) $ is the weight function for HCP.
The orthogonal condition for HCP is as follows
\begin{equation}\label{wht}
\langle T_{nm}(x),T_{n^\prime m^\prime}(x)\rangle_{\tilde{w}}=\int_{t_0}^{t_f}\tilde{w}(x) T_{nm}(x)T_{n^\prime m^\prime}(x) dx=\frac{\delta_{m,m^{\prime}}}{AN}\left\lbrace \begin{array}{lll}\vspace*{0.15in}
\pi ,~~~~n=n^\prime, m=0,\\\vspace*{0.15in}
\dfrac{\pi}{2},~~~n=n^\prime, m\neq 0,\\
0,~~~~~~~~~~~n\neq n^\prime.
\end{array}\right.
\end{equation}
\\ For brevity, let $ H_r(t):=\lbrace T_{nm}(t)\rbrace_{_{n,m}}, r=1,...,NM$.
A function $ f(t)\in L^{^2}(D)$ can be expanded  in terms of hybrid functions as
\begin{equation}\label{ff}
 f(t)\simeq f_{_{NM}}(t)=\sum\limits_{m=0}^{M-1} \sum\limits_{n=1}^{N} c_{nm} T_{nm}(t)=\sum\limits_{r=1}^{NM}c_{r} H_{r}={\bf C}^{T}{\bf H}(t)
\end{equation}
where
\[{\bf C}=[c_{_{10}},...,c_{_{1M-1}},c_{_{20}},...,c_{_{2M-1}},...,c_{_{N0}},...,c_{_{NM-1}}]^{T}=[c_{_1},c_{_2},...,c_{_{NM}}]^{T},\]
\[{\bf H}(t)=[T_{_{10}}(t),...,T_{_{1M-1}}(t),T_{_{20}}(t),...,T_{_{2M-1}}(t),...,T_{_{N0}}(t),...,T_{_{NM-1}}(t)]^{T}=[H_{_1}(t),H_{_2}(t),...,H_{_{NM}}(t)]^{T}.\]
Also, it should be noted that for $ N=1 $, HCP is equal to SCP.\\
\subsection{Function approximation}
Let $ X=L^2(D) $ and $ X_{_{NM}}=Span\lbrace H_{_1}(t),H_{_2}(t),...,H_{_{NM}}(t)\rbrace $. Since $ X_{_{NM}}\subset X$, then for every $ u \in X, $ there exist a unique best approximation of $ X_{_{NM}}$  such that
\begin{equation}\label{inf}
 \Vert u-\overline{u}_{_{NM}}(t)\Vert=\inf\limits_{g \in X_{_{NM}}} \Vert u-g\Vert,
\end{equation}
and
\begin{equation}
\overline{u}_{_{NM}}(t)=\sum\limits_{i=1}^{NM} \overline{u}_i H_i(t) =\overline{U}^T {\bf H}(t),
\end{equation}
where $ \overline{U}=[\overline{u}_1,...,\overline{u}_{_{NM}}]^{T}. $
\subsection{Operational matrix of integration and product}
\subsubsection{The SCP operational matrix of integration and product}
For convenience, let consider $  \bT(t)=[T_0(t),T_1(t),...,T_{M-1}(t)] $  as the vector of SCP basis functions for an arbitrary $ M. $
The operational matrix of integration for the shifted CP was derived by Shih (1983) [see Datta and Mohan (1995), p. 117] which satisfies in the following expression:
\begin{equation}\label{opp}
\int_{t_{0}}^{t}{\bf T}(s) ds\simeq{\bf P}~{\bf T}(t).
\end{equation}
where
\begin{equation}\label{opam}
{\bf P}=\dfrac{1}{A}\left(
\begin{array}{lllllllllr}
~~1&0 &1& 0&\cdots &0&0&0\vspace*{0.1in}\\
\frac{-1}{4}&0 &\frac{1}{4}& 0&\cdots &0&0&0\vspace*{0.1in}\\
\frac{-1}{3}&\frac{-1}{2} &~0&\frac{1}{6}&\cdots &0&0&0\\\\
\vdots & \vdots &\vdots &\vdots & &\hspace*{-0.1in}\ddots &\vdots &\vdots\vspace*{0.05in}\\
\frac{(-1)^{M-1}}{(M-1)(M-3)}&0&0&0&\cdots & \frac{-1}{2(M-3)}&0&\frac{1}{2(M-1)}\vspace*{0.05in}\\
\frac{(-1)^{^{M}}}{M(M-2)}&0&0&0&\cdots &0&\frac{-1}{2(M-2)}&0\\
\end{array}\right).
\end{equation}where $ A $ is defined by (\ref{A}).\\
The matrix  $ \bP $ for CP is obtained regarding to the following property
\begin{equation}\label{tm}
 \int_{t_{0}}^{t} T_{m-1}(s)ds= \dfrac{1}{A}(\dfrac{(-1)^m}{(m-1)^2-1}T_{0}(t)-\dfrac{1}{2(m-2)}T_{m-2}(t)+\dfrac{1}{2m}T_{m}(t)),~~~~~~~m\geq 3,
\end{equation}
where $ A=\frac{2}{t_f-t_0}.$\\\\
Operational matrix of product for CP is defined in Maleknejad et al. (2007) which is the consequences of the following property of CP,
\begin{equation*}\label{q24}
\phi_{i}(t)\phi_{j}(t)=\frac{1}{2}(\phi_{i+j}(t)+\phi_{\vert i-j \vert}(t)).
\end{equation*}
\\
Since, SCP basis functions satisfy the above relation, so SCP has the same operational matrix product.
 In fact, the product operational matrix  $ \overline{\bC}^{^T} $ satisfies in the following expression
 \begin{equation}\label{c}
 \bT(t)\bT^{^T}(t)\bC\simeq\overline{\bC}^{^T}\bT(t),
 \end{equation}
 where the vector $ \bC:=[c_{_0},c_{_1},...,c_{_{N}}] $ and the matrix $ \overline{\bC} $ is a square matrix of order N+1.\\
 \subsubsection{The HCP operational matrix of integration and product}
The vector $ {\bf H}(t)$ defined in Eq. (\ref{ff}) is a vector of hybrid Chebyshev polynomials and block-pulse functions.
The HCP operational matrix of integration over the interval $ [t_0,t_f]$ is as follows
\begin{equation}\label{op1}
\int_{t_{0}}^{t}{\bf H}(s) ds\simeq{\bf Q}~{\bf H}(t),
\end{equation}where
\begin{equation}\label{oph}
{\bf Q}=\dfrac{1}{N}\left(
\begin{array}{lllllllllr}
{\bf P}&{\bf E} &{\bf E}&\cdots &{\bf E}\\
0&{\bf P} &{\bf E}&\cdots &{\bf E}\\
0 &0 &{\bf P}&\cdots &{\bf E}\\
 \vdots &\vdots &\vdots &\ddots &\vdots\\
0&0&0&\cdots &{\bf P}\\
\end{array}\right)_{NM\times NM}
\end{equation} where $ {\bf P} $ is defined in (\ref{opam}) and
\begin{equation*}
{\bf E}=\dfrac{1}{A}\left(
\begin{array}{lllllllllr}
2& 0&0&\cdots &0\\
0&0 &0&\cdots &0\\
\frac{-2}{3} &0 &0&\cdots &0\\
 \vdots &\vdots &\vdots & &\vdots\\
\frac{(-1)^{M}-1}{M(M-2)}&0&0&\cdots &0\\
\end{array}\right)_{M\times M}.
\end{equation*}
Product operational matrix for HCP is as follows
 $$\bH(t)\bH^{^T}(t)\bC=\overline{\bC}^{^T}\bH(t)$$
 where  $ \bC=[c_{_{10}},c_{_{11}},...,c_{_{1M-1}},...,c_{_{N0}},c_{_{N1}},...,c_{_{NM-1}}]^T $ is an $ NM-$vector and
$\overline{\bC}=diag({\bf C}_{_1},{\bf C}_{_2},...,{\bf C}_{_N})$
where
 each  matrix $ {\bf C}_i $ has a similar structure to the matrix $ \bf C $ defined in (\ref{c}) and the elements of each matrix $ {\bf C}_i $ are from the vector $ [c_{_{i0}},c_{_{i1}},...,c_{_{iM-1}}] $  for $ i=1,...,N.$
\subsection{  The Product Operational Vector}\label{subb}
\vspace*{0.15in}
This section is devoted to introduce the product operational vector for HCP based on introduced operational vector of CP in (Maleknejad and Dehbozorgi 2018). For convenience, we briefly restate about this product operational vector. 
\subsubsection{The SCP operational vector }
One of the important properties of the block-pulse functions (BF) is disjointness, for more details see Jiang et al. (1992). It is yield that an n-vector of block pulse basis functions,  $  \phi=[\phi_0(t),\phi_1,...,\phi_{n-1}]$, has the following operational vector
\begin{equation*}
{\bf\phi}^T(t)~ {\bf B~\phi}(t)={\bf\hat{B} ~\phi}(t),
\end{equation*}
where the matrix $  \bf B$ is a square matrix of order n and the vector $ \bf\hat{B}$ is an n-vector with elements equal to the diagonal entries of matrix $  \bf B$.
Here,  we  introduce these operational vectors as  explicit and closed formulae for shifted Chebyshev polynomials with respect to the following property
\begin{equation}\label{q}
T_{i}(t)T_{j}(t)=\frac{1}{2}(T_{i+j}(t)+T_{\vert i-j \vert}(t)).
\end{equation}
Suppose that ${\bf B}=(b_{i,j})_{M\times M}$ is a square matrix, hence by using the above expression, we can achieve the CP operational vector  as follows:
\begin{equation}\label{bhat}
{\bf T}^{T}(t)~ {\bf B~T}(t)\simeq{\bf\hat{B} ~T}(t),
\end{equation}
where the entries of the vector ${\bf \hat{B}}$ can be interpreted as
\begin{equation}
{\bf \hat{B}}(k)=\sum \limits_{i,j=1}^{M} c_{i,j} b_{i,j},~~~~~ k=1,...,M,
\end{equation}
where
\begin{equation}\label{cij}
c_{i,j}=\left\lbrace \begin{array}{cc}\vspace*{0.15in}
1 ,~~~~\mathcal{A} \wedge \mathcal{B},\\ \vspace*{0.15in}
\frac{1}{2},~~~~~\mathcal{A}~\underline{\vee}~\mathcal{B},\\
0,~~~ otherwise.
\end{array}
\right.
\end{equation}
For brevity, two conditions $\vert i-j\vert=k-1$ and $i+j=k+1$ are defined by $\mathcal{A}$ and $\mathcal{B}$, respectively.
The symbol $\underline{\vee}$ is a logical symbol which means that  the expression  $\mathcal{A}~\underline{\vee}~\mathcal{B}$ is true if and only if just one condition is true.\\
The advantages of the prescribed vector are as follows:
\begin{itemize}
\item Interpreting the positive integer power of a function via an explicit formula as follows
\begin{equation}\label{uhat}
[u(t)] ^r\simeq {\hat{\bU}}_{r}~\bT(t),
\end{equation}
where $ u(t) \simeq u_{_N}(t):=\bU \bT(t). $\\
\item Omitting one of the basis functions vector $\bT(t)$ in such cases like (\ref{bhat}). It is a useful characteristic to transform NVIE1 to a system of algebraic  equations.
\end{itemize}
\subsubsection{The HCP operational vector}
Since HCP basis functions are included the SCP basis functions when $ N=1,$
so the SCP operational vector can be applied for constructing the HCP operational vector.
The HCP operational vector is defined as follows:
\begin{equation}\label{bhat1}
{\bf\hat{B}}=({\bf\hat{B}}_1 ,...,{\bf\hat{B}}_{N}),
\end{equation}
where each matrix $ {\bf\hat{B}}_l,~ l=1,...,N, $ is similar to the matrix $ {\bf\hat{B}} $ defined in (\ref{bhat}).
  Owing to the orthogonality of $ T_{im} $ for different $ i $, only the diagonal N-blocks of the matrix $ \bB $ are considered, i.e.
\[
~~~~~{\bf \hat{B}}_l(k)=\sum\limits_{i,j} c_{i,j} b_{i,j},~~~~~~ i,j,k=(l-1)M+1,...,lM,~~~  l=1,...,N,
\] where $ c_{i,j} $ is similarly defined as (\ref{cij}) but the conditions $ \mathcal{A} $ and $ \mathcal{B} $ are $  \vert i-j\vert=k-(l-1)M-1$ and $ i+j=k+(l-1)M+1 $, respectively.
\section{Direct method description}\label{sec2}
In this section, we illustrate how  the present scheme can be used to convert NVIE1 into a system of algebraic equations.
Note that for simplicity, our discussion and notation are restricted to Chebyshev polynomials $ \bT(t) $ but it can be easily  generalized for HCP, i.e. $ \bH(t) $.
Analogy with Maleknejad et al. (2011), we present a direct method based on operational vectors, but here we introduce an explicit formula to obtain this vector and also the proposed scheme works well for all classes of nonlinear first-kind Volterra integral equation which are inherently ill-posed.\\\\
{\bf Direct operational vector method (DOV)}\\\\
This section will introduce a direct scheme using the operational vector and matrices for approximating NVIE1. First, all functions need to be expanded with respect to CP as
$ \vspace*{-0.1in} $\[u(t)\simeq u_{_N}(t)={\bf U}^{^T}{\bf T}(t), ~~~~k(s,t)\simeq k_{_{N}}(s,t)={\bf T}^{^T}(s){\bf K}{\bf T}(t), \vspace*{-0.1in} \]
\[G(u(t)) \simeq \bZ^{^T} \bT(t).\]
Nonlinearities of function $ G $ can be classified into three types: (i) Invertible nonlinearity, (ii) Algebraic nonlinearity, (iii)  Non-invertible and non-algebraic nonlinearity. We investigate different techniques for these nonlinearities. It should be noticed that the technique which is used for the item (iii) can use for two other nonlinearity.
\subsection{Invertible nonlinearity}\label{in}
If the nonlinear function $G  $ be invertible, then
by applying the above expression for NVIE1, Eq. (\ref{V}) will be approximated as
$$\bF ^{^{T}}\bT(t)= \bT^{^{T}}(t) \int \limits_{t_{0}}^{t} \bk ~\bT(s)\bZ^{^T} \bT(s) ds,$$
 Also, by using operational matrix of integration and product, one can conclude
 \begin{equation}\label{mr}
 \bF ^{^{T}}\bT(t)= \bT^{^{T}}(t) \bk \overline{\bZ}^{^T}\bP \bT(t).
 \end{equation}
Let $ Q:=\bk \overline{\bZ}^{^T}\bP $, then our introduced operational vector simplifies the above expression as an  matrix representation
\[\bF ^{^{T}}= {\hat{Q}}.\]
Note that without the operational vector, the prevalent collocation method is also incapable to solve (\ref{mr}). Here we consider $ G $ is invertible, then from obtained $ \bZ $ and $G(u(t)) \simeq \bZ^{^T} \bT(t)$, the unknown function $ u(t) $ can be derived.
\subsection{Algebraic nonlinearity}\label{al}
For an algebraic nonlinear function $G(u(s))= \sum \limits_{r=0}^{n} \alpha_{r} u^{r}(s),$ Eq. (\ref{V}) is as follows
 \begin{equation}\label{V2}
f(t)=\sum \limits_{r=0}^{n}\alpha_{r}\int_{t_0}^t K(x,t)u^{r}(x) dx.
\end{equation}
 Using Eq. (\ref{uhat}), the  Eq. (\ref{V2}) can be approximated as
$$\bF ^{^{T}}\bT(t)= \sum \limits_{r=0}^{n} \alpha_{r} (\bT(t) \int \limits_{t_{0}}^{t} \bk ~\bT(s) ~\hat{\bU}_{r}~\bT(s) ds).$$
Due to the operational matrix of integration and product, the following result is obtained
\begin{equation}\label{eq1}
\begin{array}{cc}
\bF ^{^{T}}\bT(t)&= \sum \limits_{r=0}^{n} \alpha_{r} (\bT(t) \int \limits_{t_{0}}^{t} ~\bk ~ \overline{{\hat{\bU}} }_{r}~\bT(s) ds)\\
&= \sum \limits_{r=0}^{n} \alpha_{r} (\bT(t)  ~\bk ~ \overline{\bU }_{r}~ \bP~\bT(t)).
\end{array}
\end{equation}
Utilizing the Eq. (\ref{bhat}) with assumption $\bQ_r:=\bk \overline{\bU} _{r} \bP$ yields following
 matrix form of the Eq. (\ref{V})
\begin{equation}\label{nv}
\bF ^{^{T}}= \sum \limits_{r=0}^{n} \alpha_{r} \hat{\bQ}_{r},
\end{equation}
where $\hat{\bQ}_{r}  $ is a nonlinear vector in terms of entries of $ \bU $. Therefore, $ u_{_{N}}(t) $ can be obtained from $ u_{_{N}}(t)=\bU^{^{T}} \bT(t) $ directly.
\subsection{Non-invertible and non-algebraic nonlinearity}
\subsection*{(I) Taylor method:}
The simple structure of our proposed scheme can be conducted for other types of nonlinear continuous function $ G $. They can be approximated as $$G(u(s))\simeq \sum \limits_{r=0}^{n} \alpha_{r} u^{r}(s),$$  then  all foregoing discussions stated for algebraic nonlinearity are valid.
\subsection*{(II) Hybrid method of operational matrix and pseudospectral  collocation methods:}
First, we follow all discussion which is stated in subsection \ref{in}. After obtaining the unknown vector $ \bZ $, we use collocation points over the interval $ D $ to obtain $ u_{_N}(t) $ by solving the following nonlinear system of equations
\begin{equation}\label{ab}
G(U\bT(t_i))=\bZ^{^T} \bT(t_i),~~~~~i=1,...,n,
\end{equation}
where the appropriate collocation points $ t_i $ are chosen.
\section{Error bounds}
\begin{theorem}
(Canuto et al. 1988)
 If $u(t) \in H^{^k}_{_w}(D)$(Sobolev space) and $u_{M}(t)=\sum \limits_{r=0}^{M-1} c_{r} T_{r}(t)={\bf C}^{^{T}}{\bf T}(t)$ be the best approximation polynomials of $ u(t) $ in $L^{^2}_{_w}$-norm, then
 $$\Vert u(t)-u_M(t)\Vert_{_{L^2_w(D)}}\leq C_0 M^{-k} \Vert u(t)\Vert_{_{H^{^k}_{_w}(D)}}. $$
\end{theorem}
\begin{theorem}\label{t1}
Suppose that $ f(t)\in C^{^M}(D)  $ and $ \overline{f}_{_{NM}}=\sum\limits_{m=0}^{M-1} \sum\limits_{n=1}^{N} c_{nm} H_{nm}(t)={\bf C}^{^T}{\bf H}(t), $ where
\[{\bf C}=[c_{_{10}},...,c_{_{1M-1}},c_{20},...,c_{_{2M-1}},...,c_{_{N0}},...,c_{_{NM-1}}]^{T},\]
\[{\bf H}=[H_{_{10}},...,H_{_{1M-1}}(t),H_{20}(t),...,H_{_{2M-1}}(t),...,H_{_{N0}}(t),...,H_{_{NM-1}}(t)]^{T},\]
be the best approximate hybrid Chebyshev polynomials of $ f(t) $ in $ L^2_{\tilde{w}}(D) $, then
\begin{equation}\label{be}
\Vert f(t)-\overline{f}_{_{NM}}(t)\Vert_{L^2_{\tilde{w}}(D)}  \leq \dfrac{\gamma}{N^{M-1}M!}(\sqrt{\dfrac{\pi}{AN}}),~~~~ ~~ \gamma=\max\limits_{t\in D}\vert f^{(M)}(t)\vert.
\end{equation}
\end{theorem}
Proof. Suppose that $ f(t)=\sum\limits_{i=1}^{N}f_i(t)$ where $ f_i(t)\in C^{^M}\left[t_0+\frac{2(i-1)}{AN},t_0+\frac{2i}{AN}\right].$
Now, consider the Taylor expansion of $ f_i(t) $ as follows
\[\widehat{f}_i(t)=f_i(a_{i-1})+f^{\prime}_i(a_{i-1})(t-a_{i-1})+...+f^{^{(M-1)}}_i(a_{i-1})\dfrac{(t-a_{i-1})^{M-1}}{(M-1)!}+... ~~~t\in [a_{i-1},a_i]\]
where $ a_{i-1}:=t_0+\frac{2(i-1)}{AN} $, $ a_i:=t_0+\frac{2i}{AN} $. The truncation error of $ \widehat{f}_i(t) $ at $ Mth $ term
can be derived as
\begin{equation}\label{tr}
\vert f_i(t)-\widehat{f}_i(t)\vert \leq \vert f^{^{(M)}}_i(a_{i-1})\vert \dfrac{(t-a_{i-1})^{M}}{M!}\leq \dfrac{\gamma_i}{N^{M}M!},~~~~\gamma_i=\max\limits_{t\in[a_{i-1},a_i]}\vert f_i^{(M)}(t)\vert.
\end{equation}
Note that the last inequality of the above expression is obtained by substituting $ t=a_i.$ \\Define $ \widehat{f}(t)=\sum\limits_{i=1}^{N}\widehat{f}_i(t) $ as the Taylor expansion of $ f(t). $
Now, by regarding the concept of the best approximation hybrid functions of $ f(t) $ (\ref{inf}), triangular inequality and Eq. (\ref{tr}), we have
\begin{equation}\label{ga}
\begin{array}{lll}
\Vert f(t)-\overline{f}_{_{NM}}(t)\Vert_{_\infty} &\leq \Vert f(t)-\widehat{f}(t) \Vert_{_\infty}= \Vert \sum\limits_{i=1}^{N} (f_i(t)-\widehat{f}_i(t)) \Vert_{_\infty} \\
&\leq \sum\limits_{i=1}^{N} \Vert f_i(t)-\widehat{f}_i(t) \Vert_{_\infty} \leq \sum\limits_{i=1}^{N} \dfrac{\gamma_i}{N^{M}M!} \leq \dfrac{\gamma}{N^{M-1}M!},
\end{array}
\end{equation}
where $\gamma=\max\limits_i \gamma_i=\max\limits_{t\in[a_{i-1},a_i]}\vert f^{(M)}(t)\vert. $
\\Since $ \tilde{w}>0 $, then  the $ L^2_{\tilde{w}} $-norm of the error can be obtained by using Holder inequality and Eqs. (\ref{wht}), (\ref{ga}) as follows
\begin{equation}
\begin{array}{lll}
\Vert f(t)-\overline{f}_{_{NM}}(t)\Vert &=(\int_{t_0}^{t_f}\vert f(t)-\overline{f}_{_{NM}}\vert^2~\tilde{w}(t)~dt)^\frac{1}{2}\\
&\leq \Vert f(t)-\overline{f}_{_{NM}}(t)\Vert_{_\infty}(\int_{t_0}^{t_f} \tilde{w}(t)~dt)^\frac{1}{2}
&\leq  \dfrac{\gamma\sqrt{\frac{\pi}{AN}}}{N^{M-1}M!}.
\end{array}
\end{equation}\\
\section{Numerical examples}
In this section, the convergence behavior of solutions which are resulted from the proposed direct operational vector (DOV) method is investigated for several examples. To this end, the maximum absolute error norm  is used which is defined as follows
$$E_{\infty}=max\lbrace \vert u(t_i)-u_N(t_i)\vert , t_i \in D \rbrace.$$
For convenience, we denote the parameter $ L $ as the quantity of the used basis functions in the approximation methods.
 \begin{example}
Consider the following nonlinear Volterra integral equation of the first kind
 $$\int\limits_{0}^{t}cos(t-x)u^{\prime\prime}(x) dx=6(1-cos(x)),~~~~~~~~~u(0)=u^\prime(0)=0,$$
which has the exact solution $ u(x)=x^3.$
 \end{example}
The approximate solution using the present scheme is in high agreement with the
exact solution. The comparison of the approximate solutions by using Haar wavelet method (Singh and Kumar 2016) and our proposed method are listed in Tables 1 and 2.  It can be observed that only small size of the operational matrix is required to provide the appropriate solution. In other words, the convergence speed of our method is much more than Haar wavelet methods. For instance, our scheme has the maximum absolute error $ 1.24e-11 $ with $ L=10 $, whereas the best absolute error of Haar wavelet method has the order  $ O(10^{-7}) $ for $ L=512.$
\begin{example}\label{ex2}
As the second example, consider following equation
 $$\int\limits_{0}^{t}e^{(t-x)}ln(u(x)) dx=e^t-t-1,$$
which has the exact solution $ u(t)=e^t.$
 \end{example}
 Tables 1 and 2 demonstrate the numerical results. As we expected, our method is more accurate with less basis functions $ L$ in comparison with other previous approaches. For instance,  the absolute error of the present scheme is $ 6e-8$ for $ L=10 $,  whereas Haar wavelet method Singh and Kumar 2016) and Sinc Nystr{\"o}m  (Ma et al. 2016) algorithms obtained $ E_{_\infty}=2.1e-7,~1e-5$ with  $ L=512,~33 $ basis functions, respectively.
 \begin{example}\label{ex3}
Consider the NVIE1
 $$\int\limits_{0}^{t}e^{(t-x)}u^2(x) dx=e^{2t}-e^{t},$$
with exact solution $ u(t)=e^t.$
 \end{example}
 This problem  has been studied in  Singh and Kumar (2016), Ma et al. (2016) and Babolian and Shamloo (2008). Babolian et al. applied operational matrices of piecewise constant orthogonal functions and Laplace transform. Their obtained absolute error has the order   $\mathcal{O}(10^{-3}) $ for $ L=16. $ Tables 1, 2 and those in Babolian and Shamloo (2008) verified that the satisfactory results with  fewer basis functions ($ L $) are provided by the proposed method. It should be pointed out that although the Sinc Nystr{\"o}m method is better than the present method by increasing the value of $ L $,  it is obvious that
the convergence speed of the present method is much more than the Sinc Nystr{\"o}m method.
 \begin{example}\label{ex4}
Consider the NVIE1
 $$\int\limits_{0}^{t}(sin(t-x)+1)cos(u(x)) dx=\dfrac{tsint}{2}+sint,$$
which has the exact solution $ u(t)=t.$
 \end{example}
 The results are reported in Tables 1 and 2. These results again establish the fact that the absolute
errors of the present method is very low in comparison with the absolute errors reported in Singh and Kumar (2016) and Ma et al. (2016) with less used basis functions.
 \begin{table}[ht]
\begin{center}
\begin{small}
\caption{Absolute error of the present method for various $  L$}\vspace*{0.1in}
\begin{tabular}{ccccccccc}
\hline
\noalign{\smallskip}
$L  $&&  Ex.1&& Ex. 2&& Ex. 3 && Ex. 4  \\
\noalign{\smallskip}\hline\noalign{\smallskip}
 $2$&&$1.63e-2$&&$2.63e-1$&&$2.4e-1$   &&$5.50e-1$     \\
 $4$&&$3.29e-4$&&$2.29e-2$&&$8.32e-3$  &&$4.35e-2$   \\
$6$&&$9.71e-7$&&$1.35e-4$  &&$4.01e-5$&&$3.37e-4$  \\
$8$  &&$2.20e-10$&&$2.20e-7$&&$2.22e-7$ &&$6.54e-6$   \\
 $10$&&$1.24e-11$&&$6.24e-8$&&$1.02e-7$  &&$3.04e-8$  \\
\noalign{\smallskip}\hline
\end{tabular}
\end{small}
\end{center}
\label{ttt2}
\end{table}

 \begin{table}[ht]
\begin{center}
\begin{small}
\caption{Absolute error of the Haar wavelet method for various $  L$}\vspace*{0.1in}
\begin{tabular}{ccccccccc}
\hline\noalign{\smallskip}
$L  $&&  Ex.1&& Ex. 2&& Ex. 3 && Ex. 4  \\
\noalign{\smallskip}\hline\noalign{\smallskip}
 $4$&&$8.4e-3$&&$3.4e-3$&&$2.8e-3$   &&$1.2e-3$     \\
 $8$&&$2.1e-3$&&$8.4e-4$&&$7.3e-4$  &&$3.1e-4$   \\
$16$&&$5.5e-4$&&$2.1e-4$  &&$1.8e-4$&&$8.0e-5$  \\
$32$  &&$1.4e-4$&&$5.5e-5$&&$4.6e-5$ &&$2.0e-5$   \\
 $64$&&$3.5e-5$&&$1.3e-5$&&$1.1e-5$  &&$5.0e-6$  \\
$128$&&$8.8e-6$&&$3.4e-6$  &&$2.9e-6$&&$1.2e-6$  \\
$256$  &&$2.2e-6$&&$8.6e-7$&&$7.3e-7$ &&$3.1e-7$   \\
 $512$&&$5.5e-7$&&$2.1e-7$&&$1.8e-7$  &&$7.9e-8$  \\
\noalign{\smallskip}\hline
\end{tabular}
\end{small}
\end{center}
\label{ttt3}
\end{table}
\begin{example}\label{ex5}
As an another test problem, consider the following NVIE1
 \begin{equation}
f(t) =\int_0^t u^3(x)dx,
 \end{equation}
 whose $f(t)$ is determined by noting that the non-smooth solution $ u(t)=\vert t-\frac{1}{2}\vert $,
 hence $f(t)$ is  $\frac{1}{64}+\frac{1}{4}(t-\frac{1}{2})^3\vert t-\frac{1}{2}\vert$ .\end{example}
Since  the unknown function $ u(t)\in C^{^0}\setminus C^{^1}$, then hybrid functions (HCP)  have a better efficiency rather than CP. When $M=4,~ N=2$ or $ L=8 $, HCP achieve $ E_{_\infty}=2.51e-14 $ vs. the absolute error of Chebyshev polynomials with $ L=8 $ have $ E_{_\infty}=7e-2$. Figure \ref{fi3} depicts the comparison of the exact and approximate solutions with the basis functions CP and  HCP. Note that in this example, $ G(u(t))=\vert t-\frac{1}{2}\vert^3 $ is a non-invertible function when $ t \in [0,1]. $ Thus, it may solve this problem by using (\ref{nv}) or (\ref{ab}).
\begin{figure}
\begin{center}
\subfigure[The HCP approximate solution and the exact solution.]{
\resizebox*{5.7cm}{!}{\includegraphics{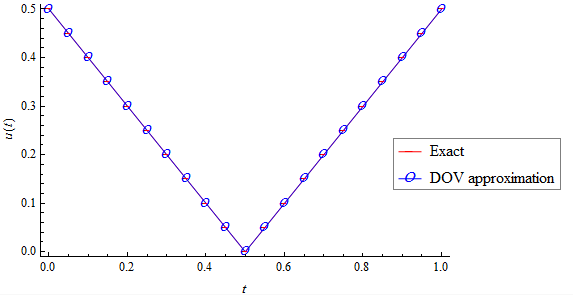}}}\hspace{3pt}
\subfigure[The CP approximate solution and the exact solution.]{
\resizebox*{5.7cm}{!}{\includegraphics{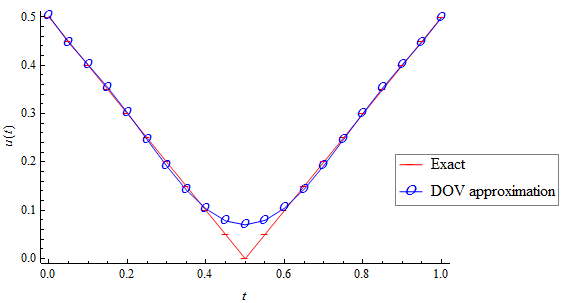}}}
\caption{Results of Ex. \ref{ex4} for L=8.}
\label{fi3}
\end{center}
\end{figure}
 \begin{example}
Consider the following linear Volterra integral equation of the first kind
 $$\int\limits_{0}^{t}e^{(t+x)}u(x) dx=te^{t},~~~~~~~~~~~t\in[0,1]$$
with exact solution $ u(t)=e^{-t}.$
 \end{example}
 This linear example has been considered in some previous works, Babolian and Masouri (2008), Masouri et al. (2010), Maleknejad et al. (2011) and Khan et al. (2014). Table 3 represents the superiority of our method respect to other studies.
  \begin{table}[ht]
\begin{center}
\begin{small}
\caption{Comparison of the absolute error of some recent methods and the presented DOV method for Ex. 6}\vspace*{0.1in}
\begin{tabular}{cccccc}
\hline
\noalign{\smallskip}
& Direct method & EI method   & BP method   & OHAM & DOV\\
&Babolian and Masouri (2008) & Masouri et al. (2010) &Maleknejad et al. (2011) & Khan et al. (2014)&\\
&($ L = 64 $) &($ L = 64 $)&($ L = 64 $) & (order 5)&($ L = 8 $)\\
\noalign{\smallskip}\hline\noalign{\smallskip}
$ E_{\infty} $&$ 1.0e-3 $&$ 1.9e-4 $ &$ 5.5e-3 $&$ 4.7e-6 $ &$ 1.29e-8 $\\
\noalign{\smallskip}\hline
\end{tabular}
\end{small}
\end{center}
\label{ttt4}
\end{table}
\begin{example}
As a test problem, consider the following NVIE1
 $$\int\limits_{0}^{t}(u^2(x)-u(x)) dx=\dfrac{t^3}{3}-\dfrac{t^2}{2}, ~~~~ t \in [0,2]$$
which has the exact solution $ u(t)=t.$
 \end{example}
In this example, the nonlinear part $ G(u(t)) $  is not invertible for all $ t\in[0,2] $. The  present scheme provides the exact solution by using only $ L=3 $ basis functions. It verifies high convergence rate and low computational complexity of the scheme.
\begin{example}\label{ex9}
As the final test problem, consider the following NVIE1
 $$\int\limits_{-1}^{t}t~x ~u^2(x) dx=f(t), ~~~~ t \in [-1,1]$$
where $ f(t) $ can be determined such that
\[u(t)=\left\lbrace\begin{array}{lll}
\sqrt{\vert t\vert},~~~~-1\leq t\leq0\vspace*{0.1in}\\
t(1-t),~~~~0\leq t\leq1
\end{array}\right.\]
 \end{example}
Fig. 2 (a) depicts the exact solution and the approximate solution $ u_{2,8} $ of this example. It is observed that this piecewise function is non-invertible. We apply the methods (\ref{nv}) or (\ref{ab}) to obtain the best approximate function. Furthermore, hybrid functions work well for these kinds of piecewise functions belongs to $ C^{0}\setminus C^{1} $. The present approach gives an approximate solution of order $ \mathcal{O}(10^{-12})$ when $ L=16. $ Fig. 2 (b) shows the approximation error of $  u_{2,8} .$\\\\
\begin{figure}
\begin{center}
\subfigure[The HCP approximate solution $  u_{2,8} $ and the exact solution.]{
\resizebox*{5.7cm}{!}{\includegraphics{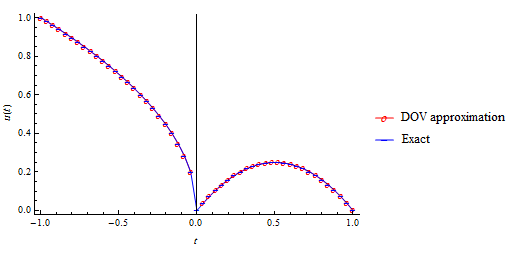}}}\hspace{3pt}
\subfigure[Comparison of the error function $  e_{2,8}=u(t)-u_{2,8}.  $ ]{
\resizebox*{5.7cm}{!}{\includegraphics{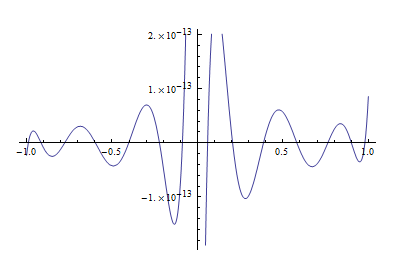}}}
\caption{Results of Ex. \ref{ex9} for $  L=16$.}
\label{fi3}
\end{center}
\end{figure}
\begin{example}
As the test problem, consider the following NVIE1
 $$\int\limits_{1}^{t}e^{u(x)} dx=cos(1)-cos(t), ~~~~ t \in [1,2]$$
which has the exact solution $ u(t)=Ln ~sin(t).$
 \end{example}
 In this example, the nonlinear function $ G(u(t))=sin(t) $ over the interval $ [1,2]$ which is a non-invertible and non-algebraic nonlinear function in this interval. Hence, we apply (\ref{ab}). Table 4 demonstrates the results for different orders of $ N $ and $ M$ which are adjustable to obtain the best approximate solution. The quantity of the basis functions is $ L=MN. $
  \begin{table}[ht]
\begin{center}
\begin{small}
\caption{Comparison of the absolute error of some recent methods and the presented DOV method for Ex. 8. }\vspace*{0.1in}
\begin{tabular}{ccccccccc}
$ N=1 $\\
\hline
\noalign{\smallskip}
 M & $ 2 $ & $4 $  &$ 6 $   &$8 $ \\
\noalign{\smallskip}\hline\noalign{\smallskip}
$ 1.0 $&$ 8.28e-04 $&$ 2.19e-06$ &$ 2.92e-09 $&$ 3.20e-09 $ \\
$ 1.2 $&$ 4.24e-04 $&$ 5.05e-07 $ &$ 1.11e-9 $&$ 1.65e-09 $ \\
$ 1.4 $&$ 1.48e-04 $&$ 1.00e-07 $ &$7.90e-10 $&$1.30e-09 $ \\
$1.6$&$ 4.17e-04 $&$ 9.90e-07 $ &$ 9.46e-10 $&$ 3.91e-10 $ \\
$ 1.8$&$1.80e-04 $&$ 1.44e-07 $ &$ 8.17e-10 $&$ 2.3e-09 $ \\
$2.0 $&$ 6.22e-03 $&$ 4.58e-06 $ &$ 1.03e-08 $&$ 1.51e-09 $ \\
\hline
\noalign{\smallskip}
$E_{_\infty} $&$ 8.28e-04 $&$ 4.58e-06 $ &$ 1.03e-08 $&$ 2.02e-09$ \\
\hline
\\
\\
$ N=2 $\\
\hline
\noalign{\smallskip}
 M & $ 2 $ & $4 $  &$ 6 $   &$8 $\\
\noalign{\smallskip}\hline\noalign{\smallskip}
$ 1.0 $&$ 2.09e-04 $&$ 1.63e-07 $ &$ 5.86e-11 $&$ 1.14e-11 $ \\
$ 1.2 $&$ 8.51e-07 $&$ 3.72e-08 $ &$ 2.60e-11 $&$5.56e-13 $ \\
$1.4$&$ 2.04e-05 $&$ 1.54e-07 $ &$ 3.94e-11 $&$3.09e-12 $ \\
$ 1.6$&$ 1.10e-04 $&$ 1.69e-07 $ &$ 4.77e-11 $&$4.33e-12 $ \\
$ 1.8 $&$ 1.97e-04 $&$1.33e-07 $ &$ 4.94e-11 $&$ 6.36e-12 $ \\
$ 2.0 $&$ 5.98e-04 $&$ 1.10e-06 $ &$ 8.90e-11 $&$1.85e-11 $ \\
\hline
\noalign{\smallskip}
$E_{_\infty} $&$ 5.98e-04 $&$1.10e-06 $ &$ 1.03e-10 $&$ 1.85e-11$ \\
\hline
\end{tabular}
\end{small}
\end{center}
\label{ttt4}
\end{table}

\begin{example}\label{ex10}
As the final test problem, consider the following linear VIE1 with discontinuous solution
 $$\int\limits_{-\frac{1}{2}}^{t}t ~x~\sqrt{u(x)} dx=f(t), ~~~~ t \in [-\frac{1}{2},1]$$
whose $ f(t) $ can be derived so that the exact solution
\[u(t)=\left\lbrace\begin{array}{lll}
e^{-2t},~~~~~~~-\frac{1}{2}\leq t<0,\vspace*{0.1in}\\
t^2,~~~~~~~~~~~0\leq t<\frac{1}{2},\vspace*{0.1in}\\
\dfrac{1}{t},~~~~~~~~~~~\frac{1}{2}\leq t\leq1.\\
\end{array}\right.\]
 \end{example}
 Hybrid functions allow us to approximate discontinuous solutions as well. In real word problems, the solution $ u(t) $ are almost piecewise functions. Moreover, if $ k(x,t)$ is a continuous function in Eq. (\ref{V}), then the behavior of functions $ f(t) $ and $ u(t) $ are almost the same. Hence, we propose hybrid functions when the function $ f(t) $ is a piecewise function and $ k(x,t) $ is continuous. In this test problem, various values of $ M $ and $ N $ are tested and the best computational result is obtained with $ u_{3,12} $ which has the absolute error of order $ \mathcal{O}(10^{-11}) $. Fig. 3 shows the efficiency of our scheme for discontinuous functions. Table \ref{tttt} indicates that the error is significantly decreased when the polynomials' degree of each subinterval increase.
   \begin{table}[ht]
\begin{center}
\begin{small}
\caption{The absolute error of different degree $ M $ when $ N=3 $ for Ex. \ref{ex10}. }\vspace*{0.1in}
\begin{tabular}{ccccccccc}
$ N=3 $\\
\hline
\noalign{\smallskip}
 M  & $4 $  &$ 6 $   & $ 8 $ &$10$ &$12 $ \\
\noalign{\smallskip}\hline\noalign{\smallskip}
$ -0.50 $&$9.74e-04 $&$ 1.42e-06$ &$8.60e-10 $&$ 1.35e-13 $ &$ 7.28e-14 $\\
$ -0.25 $&$ 1.62e-05 $&$ 8.97e-09 $ &$ 2.57e-12 $&$3.33e-15 $&$ 3.99e-15 $ \\
$ 0.00 $&$ 9.37e-02 $&$2.74e-04 $ &$2.75e-07 $&$1.45e-09 $ &$1.01e-09$\\
$0.25$&$ 1.10e-03 $&$ 1.18e-06 $ &$ 5.62e-10 $&$2.08e-12 $ &$ 9.16e-13 $\\
$ 0.50$&$2.42e-03 $&$ 6.06e-07 $ &$ 8.26e-08 $&$2.25e-9 $ &$ 6.60e-11 $\\
$0.75 $&$ 3.54e-04 $&$1.76e-07 $ &$ 4.12e-09 $&$9.29e-11 $&$ 2.47e-12 $ \\
$ 1.00 $&$ 1.46e-03 $&$ 2.02e-06 $ &$1.59e-09 $&$6.54e-11 $&$ 4.77e-12 $ \\
\hline
\noalign{\smallskip}
$E_{_\infty} $&$ 9.37e-02 $&$ 2.74e-04 $ &$ 2.75e-07 $&$ 2.25e-9$&$ 9.74e-11 $ \\
\hline
\end{tabular}
\end{small}
\end{center}
\label{tttt}
\end{table}
 \begin{figure}
\begin{center}
\subfigure[The exact solution of Ex. \ref{ex10}.]{
\resizebox*{5.7cm}{!}{\includegraphics{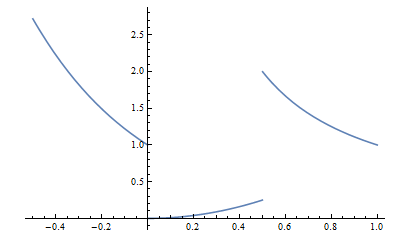}}}\hspace{3pt}
\subfigure[The error function $  e_{3,12}=u(t)-u_{3,12}(t).  $ ]{
\resizebox*{5.7cm}{!}{\includegraphics{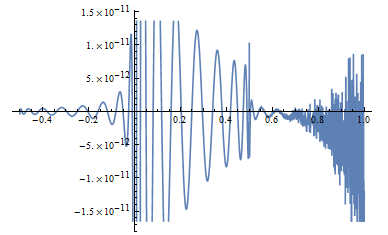}}}
\caption{Results of Ex. \ref{ex10} for $ u_{3,12}$.}
\label{fi3}
\end{center}
\end{figure}
\\
\\ In all above examples, it can be observed that for small $ L $, the approximate solutions with high accuracy are obtained. Consequently, in analogy with other methods represented in Masouri et al. 2010, Babolan and Masouri 2008, Maleknejad et al. 2011, Singh and Kumar 2016, Ma et al. 2016, we require solving  a small nonlinear system of algebraic equations to obtain an appropriate solution.
\paragraph{Remark.}It is noticeable that the invertibility of $ G(u(t)) $ is a main constraint in the previous works (Singh and Kumar 2016, Ma et al. 2016). This limitation on $ G(u) $ may be omitted using the present method. For instance in Ex. 5, 7, 8 and 9, $ G(u) $ is non-invertible.
\section{Conclusions}\label{con1}
  Explicit formulas for operational vectors have been derived based on Chebyshev polynomials. 
These vectors allow us to introduce an efficient, accurate and reliable numerical method which works well for nonlinear Volterra integral equation of the first kind.  The problem has been reduced to solving a set of algebraic equations. The main advantages of this method are ease of comprehending, simplicity of performing, high accuracy and appropriate convergence rate. In comparison with other numerical schemes such as proposed in Singh and Kumar (2016), Ma et al. (2016), the main properties of our proposed method are low storage
requirement and computational complexity with high precision of the suggested procedure. Furthermore, hybrid functions allow us to adjust the order of polynomials' degree and block-pulse functions to achieve the best computational results, especially when the unknown solution belongs to $ C^{0}\setminus C^{1} $ or discontinuous functions.
 Numerical experiments confirm that our proposed method is a simple and a powerful tool to conquer the ill-posedness and the nonlinearity of these problems.

\end{document}